\documentclass[smallextended,natbib 
]{svjour3}

\smartqed

\usepackage{graphicx,color}
\usepackage{times}

\usepackage{lineno}

\usepackage[colorlinks=true, breaklinks=true, pdfstartview=FitV, linkcolor=red, citecolor=blue, urlcolor=black]
{hyperref}

\usepackage[T1]{fontenc}

\newcommand{\Vec}[1]{\mbox{\boldmath$#1$}}

\journalname{Celestial Mechanics and Dynamical Astronomy}

\begin{document}

\title{ A new radial, natural, higher order intermediary of the main problem four decades after the elimination of the parallax. }

\author{Martin Lara}

\institute{M.~Lara \at
           GRUCACI, University of La Rioja, and Space Dynamics Group -- UPM \\
           \email{mlara0@gmail.com}
}

\date{This is a post-peer-review, pre-copyedit version of an article published in Celestial Mechanics and Dynamical Astronomy. The final authenticated version is available online at: \\ \href{http://dx.doi.org/10.1007/s10569-019-9921-5}{http://dx.doi.org/10.1007/s10569-019-9921-5}}

\maketitle

\begin{abstract}
Simplifications in dealing with the equation of the center when the short-period effects are removed from a zonal Hamiltonian are commonly attributed to the elimination of parallactic terms. But this interpretation is incorrect, and the simplifications rather stem from the removal of concomitant long-period terms, an outcome that can also be achieved without need of eliminating the parallax. To show that, a Lie transforms simplification is invented that augments the exponents of the inverse of the radius and still achieves analogous simplifications in handling the equation of the center to those provided by the classical elimination of the parallax simplification. The particular case in which the new transformation does not modify the exponents of the parallactic terms of the original problem, leads to a new intermediary of the main problem that, while keeping higher order effects of $J_2$, is formally analogous to Cid's first order radial intermediary.
\end{abstract}

\keywords{perturbed Keplerian motion \and Lie transforms \and Hamiltonian simplification \and elimination of the parallax \and elimination of the perigee \and radial intermediaries \and main problem }


\section{Introduction}

The classical way of approaching perturbed Keplerian motion is based on the explicit appearance of the mean anomaly in the disturbing function \citep{Hansen1855,Delaunay1860,Tisserand1889,BrouwerClemence1961}. The perturbation is thus represented in the form of a multivariate Fourier series in the action-angle variables of the Kepler problem, the arguments of which are linear combinations of the angles. However, a well known drawback of this approach in artificial satellite theory, is that it may provide the short-period elements of the solution with unwieldy series \citep{DepritRom1970,Kinoshita1977,Claes1980}.\footnote{The impossibility of integrating the equation of the center in closed form within the algebra of trigonometric functions \citep{Jefferys1971} caused, for sometime, the belief that expansions in powers of the eccentricity were unavoidable in a higher order normalization of the zonal problem. Still, the later invention of Hamiltonian simplification techniques allowed to relegate this issue to higher orders of the short-period elimination \citep{Deprit1981,DepritJGCD1981,CoffeyDeprit1982,AlfriendCoffey1984}. On the other hand, an early comment by \citeauthor{Aksnes1971} (\citeyear{Aksnes1971}) motivated subsequent research that showed that direct integration of the equation of the center may be, in fact, unneeded. Indeed, it can be avoided by finding the antiderivatives of the proper combination of the equation of the center with other periodic functions that appear in the short-period elimination \citep{Deprit1982,Healy2000}. Alternatively, the equation of the center can be integrated alone if the algebra of the short-period elimination is extended to encompass special functions of the polylogarithmic type \citep{OsacarPalacian1994}.} Partial remedy to this inconvenience is found by avoiding expansions in powers of the eccentricity, on the one hand, and by grouping these corrections using the known relations between the elements of the ellipse and the polar coordinates, on the other \citep{Brouwer1959}, or, much better, by transforming the short-period corrections to the orbital elements, once they have been computed, by equivalent, but notably more compact short-period corrections to the polar coordinates \citep{Kozai1962}.
\par

But if at the end, \cite{Izsak1963AJ} asked, the short-period corrections are to be given in polar coordinates, why the need of recomputing them from corresponding corrections in elements instead of deriving these corrections directly from the generating function? Izsak's point of view was supported by the fact that the polar coordinates jointly with the radial velocity, the total angular momentum, the right ascension of the ascending node, and the third component of the angular momentum vector comprise a set of canonical variables that had been successfully used by \cite{Hill1913} in the treatment of the planetary problem \cite[see also][]{Aksnes1972}. When dealing with the zonal problem of artificial satellite theory \citep{Brouwer1959,Kozai1959}, the task is further simplified by Izsak's observation that Brouwer's generating function of the short-period elimination is easily rewritten in polar variables except for the equation of the center, whose required partial derivatives with respect to the polar variables were, nevertheless, obtained by Izsak from usual expressions of Keplerian motion.\footnote{An elegant relation between the true and eccentric anomalies was some years later derived that can be used to express the equation of the center as a function of the true anomaly alone \citep{BrouckeCefola1973}, in this way making the formulation of the equation of the center in polar variables trivial.} 
\par

A further step in the same direction was given by \cite{CidLahulla1969}, who took the point of view of normalizing the perturbation problem directly in polar-nodal variables by removing the argument of the latitude from the Hamiltonian, as opposed to the traditional averaging over the mean anomaly. {This way, instead of the semi-major axis, it is the total angular momentum that becomes the formal integral of the averaged problem}. When applied to the main problem of the artificial satellite theory, the new procedure is definitely successful at the first order of the perturbation approach, in which the oblateness coefficient is taken as the small parameter. Indeed, the short-period corrections adopt a simpler form than in the traditional treatment ---a result that is, in fact, expected because the mean anomaly is removed only partially from the normalized Hamiltonian, in which some short-period terms still remain. Besides, the resulting, integrable Hamiltonian keeps most of the essence of the main problem because it is formally equivalent to a simple truncation of the Hamiltonian in which the removed terms have no direct effect in the secular and long-period dynamics \citep{DepritFerrer1987}. Initial efforts in computing the short-period corrections to the second order were, nevertheless, less appealing because mixed secular-periodic terms arose in the computations \citep{CidLahulla1971}. However, the undesired appearance of these non-periodic terms is not related to the procedure itself, {but to the particular} choice by Cid and Lahulla of the null value for the arbitrary ``constant'' that appears in the computation of the scalar generating function from which the periodic corrections are derived.
\par

The indetermination of the generating function of the perturbation approach up to an arbitrary function, which in the case of perturbed Keplerian motion only needs to meet the condition of having null derivative with respect to the mean anomaly, had been identified as early as in \citep{Kozai1962}. Later, these kinds of arbitrary functions were routinely used in the computation of centered mean elements to better represent the average evolution of the osculating elements \citep[see][for instance]{MetrisExertier1995,Steichen1998,LaraSanJuanFolcikCefola2011}. However, it belongs to \cite{AlfriendCoffey1984} the merit of demonstrating the important role that these integration ``constants'' play when the Hamiltonian reduction is attacked in other variables than the customary action and angles.
\par

Indeed, in the general framework provided by Lie series and the Lie transforms method \citep{Grobner1960,Hori1966,Deprit1969,HenrardRoels1974,GiorgilliGalgani1978,GiorgilliGalgani1985}, Alfriend and Coffey take advantage of a remarkable simplification technique devised by \cite{Deprit1981} ---the elimination of the parallax--- and propose a diversion from the usual perturbation approach. Instead of removing the remaining short-period terms from the simplified Hamiltonian obtained after the elimination of the parallax, they remove the long-period terms, which are related to the perigee dynamics, previously to approaching the complete elimination of short-period effects by the standard Delaunay normalization \citep{Deprit1982}. In this process, the proper selection of the arbitrary function that appears at each order of the Lie transforms procedure in the solution of the homological equation, is crucial to preventing the appearance of secular terms in the generating function of \citeauthor{AlfriendCoffey1984}'s (\citeyear{AlfriendCoffey1984}) \emph{elimination of the perigee} transformation, and delays by one order the correct computation of the generating function with respect to the order in which the Hamiltonian term with the argument of the perigee removed is determined. Alfriend and Coffey's approach of eliminating consecutively the parallax and the perigee, and, in consequence, achieving the removal of the argument of the latitude from the zonal Hamiltonian, was useful also in demonstrating that the elimination of the latitude is not unique \citep{CidFerrerSeinEchaluce1986}. 
\par

On the other hand, in view of the group structure of canonical transformations \citep[see][for instance]{Arnold1989,GoldsteinPooleSafko2001}, it seems natural to pursue the elimination of the latitude by a single Lie transformation as opposed to Alfriend and Coffey's decomposition. Theoretical efforts by \cite{DepritFerrer1989} to furnish the elimination of the latitude simplification with algebraic foundations created some controversy due to the variety of radial intermediaries that can be obtained in the process of eliminating the latitude \citep{CoppolaPalacian1994}. Besides, from the viewpoint of a practitioner, the straightforward composition of the sequence formed by the elimination of the parallax and the consequent elimination of the perigee into a single Lie transformation produces a notable increase in the size of the transformation equations, with the consequent deterioration in the evaluation of the periodic terms of the solution, cf.~\citep{SanJuanOrtigosaLopezOchoaLopez2013}. This latter result was in fact expected due to the clear advantages provided by the \textit{divide et impera} strategy, which is central to Deprit's Hamiltonian simplification concept.
\par

All these efforts in developing specific simplification algorithms that expedite the computation of the periodic corrections of the different solutions to the artificial satellite problem, followed the mainstream originated with the elimination of the parallax. That is, the computation of the perturbation solution by a Lie transformation is approached in polar-nodal variables within the particular algebra of the state functions defined by the conic parameter and the projections of the eccentricity vector in the nodal frame. However, it has been the sustained view of the author along the last years that the direct formulation in Delaunay canonical variables, rather than in the polar-nodal set, strips the simplification procedures of unnecessary intricacies and makes notably easier the implementation of the different simplification algorithms \citep{LaraSanJuanLopezOchoa2013,LaraSanJuanLopezOchoa2013c,LaraSanJuanLopezOchoa2013b}. And it seems to be exactly the case with the elimination of the latitude simplification. Actually, the reported multiple representations of functions in the elimination of the latitude is just a consequence of the freedom in the selection of the integer exponent with which the radius can be predetermined to stay in the simplified Hamiltonian, a point that becomes obvious when the radial distance is treated as a function of the canonical variables rather than a variable by itself. Note that approaching the Hamiltonian simplification  in Delaunay variables does not deny the usefulness of the polar-nodal variables in the representation of the periodic corrections, which, as previously mentioned, provide real advantages when untied from the state functions of the polar variables representation. 
\par

This latter fact ---the freedom in choosing the exponent of the inverse of the radius in the simplified zonal potential, without constraining to the case of the removal of parallactic factors--- provides the motivation for looking back to Deprit's elimination of the parallax. Further than the role it plays in the elimination of the latitude, it is worth checking which is the real origin of the simplification it provides in dealing with the essential short-period terms that are associated with the equation of the center.
\par

To do that, the main problem of artificial satellite theory is chosen as a demonstration model for simplicity, but extending the procedure to a full zonal potential will not introduce other additional complications than those derived from the manipulation of longer trigonometric series. Indeed, zonal perturbations just comprise parallactic factors and Legendre polynomials, thus accepting the same axial symmetry of the main problem. Hence their algebraic structure is amenable to the same treatment applied to the simpler truncation of the zonal potential to the second degree \citep{Brouwer1959,AlfriendCoffey1984,CoffeyDepritDeprit1994}. The usual formulation of the main problem in orbital elements, as well as the basics of the Lie transforms method, as reported by \cite{Deprit1969}, are recalled for completeness in Section \ref{s:preliminaries}. Then, after illustrating in Section \ref{s:Brouwer} the difficulties introduced by Brouwer's direct approach of averaging the Hamiltonian over the mean anomaly, the simplifying effects resulting from the elimination of the parallax in the elimination of short-period effects are reviewed at the beginning of Section \ref{s:simcom}. In this Section a new Lie transforms simplification, based on the elimination of the parallax, is presented that increases the exponents of the inverse of the radius, contrary to the elimination of the parallactic terms of the disturbing function by reducing them. Because the new transformation also eliminates the latitude at the first order, it is found that it provides analogous simplifications to those supplied by the elimination of the parallax in dealing with the equation of the center, yet it leads to longer series in standard trigonometric terms that handicap the short-period corrections. The fact that the new transformation does not reduce the exponents of the parallactic factors, shows that the simplifications in dealing with the equation of the center, which are achieved with either the elimination of the parallax or the new Lie transformation, are of a general nature. It will be shown that the simplification is associated to the decoupling of short-period effects of the disturbing function in the out-of-plane direction.
\par

While the Lie transformation used in Section \ref{s:simcom} has been chosen just for illustration purposes, and it has not been found to provide any advantage over the elimination of the parallax, a different case of the same approach is presented in Section \ref{s:Cid} that keeps the exponents of the radius as in the original main problem. It is found that this variant provides shorter periodic corrections at the first order than those obtained in the same order of the elimination of the parallax, while the second order corrections remain of analogous complexity in both canonical transformations. Besides, if this new transformation is complemented with Alfriend and Coffey's elimination of the perigee, in which the same constraint in what respects to the exponent of the radius is used, a higher order intermediary solution is obtained that keeps more features of the original main problem than a simpler quasi-Keplerian system. In particular, after neglecting terms of the order of the square of the eccentricity and higher, the new higher order intermediary can be cast in the simple form of Cid's first order original intermediary, whose solution is standard \citep{Ferrandiz1986,LaraGurfil2012}.
\par

Finally, it is worth mentioning that the approach taken here of using Delaunay variables in the computation of the elimination of the latitude is not in contradiction with the economy produced in the short-period corrections. Indeed, while this approach seems to depart from Cid and Lahulla's original plan, the transformations can still be derived from the generating function in polar variables in the way proposed by Izsak. In addition, if the latitude is eliminated by preserving the third power of the inverse of the radius, as reported in Section \ref{s:Cid}, the corrections to the polar variables provide additional savings with respect to existing corrections in the literature \citep{Lara2015ASR}. 
\par

To conclude this Introduction, it is mandatory to remark that the method of Lie transformations and normal form theory are not constrained at all to the old-fashioned topic of finding solutions of the main problem of artificial satellite theory. Quite on the contrary, the strength of this methodology have been applied with great success in different fields of celestial mechanics, as it is demonstrated by the wealth of contributions on the topic, among which mention could be done to the works of \citep{SansotteraCeccaroni2017,HautesserresLaraCeMDA2017,GachetCellettiPucaccoEfthymiopoulos2017,PaezEfthymiopoulos2018,MahajanVadaliAlfriend2018}, just to mention a few recent examples.
\par

\section{Dynamical model and perturbation method} \label{s:preliminaries}

The main problem of artificial satellite theory \citep{Brouwer1959} is the most simple, nonintegrable truncation of the geopotential that, in spite of it does not take into account a realistic mass distribution of the Earth, provides a useful approximation of the full problem. When written in Hamiltonian form, the main problem is \citep{BrouwerClemence1961}
\begin{equation} \label{Hammain}
\mathcal{H}=-\frac{\mu}{2a}+\frac{\mu}{r}\frac{R_\oplus^2}{r^2}\frac{1}{2}C_{2,0}\left(1-\frac{3}{2}s^2+\frac{3}{2}s^2\cos2\theta\right),
\end{equation}
where $\mu$ is the Earth's gravitational parameter, $R_\oplus$ its mean equatorial radius, $C_{2,0}=-J_2$ is the zonal harmonic coefficient of the second degree, and the orbit semi-major axis $a$, the radius $r$, the sine of the inclination $s\equiv\sin{I}$, and the argument of the latitude $\theta$ are functions of some set $(\Vec{x},\Vec{X})$ of canonical variables. In particular, when using polar variables $r$ and $\theta$ are coordinates by themselves, whose conjugate momenta are the radial velocity $R$ and the total angular momentum $\Theta$, respectively, whereas the Keplerian part, given by the first summand of Eq.~(\ref{Hammain}), is
\begin{equation} \label{keplerian}
-\frac{\mu}{2a}=\frac{1}{2}\left(R^2+\frac{\Theta^2}{r^2}\right)-\frac{\mu}{r},
\end{equation}
and the sine of the inclination is
\[
s=\sqrt{1-\cos^2I}=\sqrt{1-N^2/\Theta^2},
\]
where $N$ is the third component of the angular momentum vector, which is an integral of the main problem because its conjugate variable $\nu$, the right ascension of the ascending node, is ignorable in Hamiltonian (\ref{Hammain}).

If, on the contrary, Delaunay canonical variables are used, then $a=L^2/\mu$, where the Delaunay action $L$ is the conjugate momentum to the mean anomaly $\ell$. The radius $r$ is given by the conic equation
\begin{equation} \label{conic}
r=\frac{p}{1+e\cos{f}},
\end{equation}
where the orbit parameter is $p=G^2/\mu$ and the total angular momentum $G=\Theta$ is the conjugate momentum to the argument of the perigee $g$, the eccentricity $e=\sqrt{1-\eta^2}$ is written in terms of the eccentricity function $\eta=G/L$, the true anomaly $f$ is an implicit function of the mean anomaly, whose computation requires inversion of the Kepler equation, $\theta=f+\omega$, with $\omega=g$, and the inclination is $I=\arccos(H/G)$, where now $H$ notes the third component of the angular momentum vector, and its conjugate (cyclic) variable, the right ascension of the ascending node, is denoted $h$.
\par

Approximate solutions to the 2 degrees of freedom Hamiltonian (\ref{Hammain}) are customarily computed by perturbation methods. In particular, in the framework provided by the Lie transforms method \citep{Deprit1969}, the main problem Hamiltonian is cast in the form of a Taylor series
\begin{equation} \label{Hold}
\mathcal{H}=\sum_{m\ge0}\frac{\epsilon^m}{m!}H_{m,0}(\Vec{x},\Vec{X}),
\end{equation}
where the zeroth order is the Keplerian term, the disturbing function $H_{1,0}$ is given by the second summand of Eq.~(\ref{Hammain}), $H_{m,0}=0$ for $m>1$, and the small parameter $\epsilon$ is taken here as formal ($\epsilon=1$), although taking $\epsilon=J_2$ is also customary. Then, the reduction of the perturbation Hamiltonian (\ref{Hold}) to a simpler or integrable expression is found by applying to it a canonical transformation
\begin{equation} \label{Liet}
\mathcal{T}:(\Vec{x},\Vec{X})\longrightarrow(\Vec{x}',\Vec{X}';\epsilon),
\end{equation}
which is derived from a generating function
\begin{equation} \label{Gen}
\mathcal{W}=\sum_{m\ge0}\frac{\epsilon^m}{m!}W_{m+1}(\Vec{x},\Vec{X}).
\end{equation}
When the transformation (\ref{Liet}) is applied to the original Hamiltonian (\ref{Hold}), it is converted into the new Hamiltonian
\begin{equation} \label{Hnew}
\mathcal{K}=\mathcal{H}(\Vec{x}(\Vec{x}',\Vec{X}';\epsilon),\Vec{X}(\Vec{x}',\Vec{X}';\epsilon);\epsilon)\equiv\sum_{m\ge0}\frac{\epsilon^m}{m!}H_{0,m}(\Vec{x}',\Vec{X}').
\end{equation}
The convenience of using a double subindex notation for the Hamiltonian terms is related to Deprit's devising of the Lie transforms algorithm, which is summarized with the recurrence formula in Eq.~(\ref{triangle}).
\par

The way in which the generating function (\ref{Gen}) is constructed depends on the simplification requirements imposed on the perturbation solution ---typically, the removal of periodic effects. This construction is done stepwise by solving the homological equation
\begin{equation} \label{homoeq}
\mathcal{L}_0(W_m)=\widetilde{H}_{0,m}-H_{0,m},
\end{equation}
in which  the Lie derivative $\mathcal{L}_0$ of the generating function is
\begin{equation} \label{LieDerivative}
\mathcal{L}_0(W_m)=\{\mathcal{H}_{0,0};W_m\},
\end{equation}
and the curly brackets stand for the Poisson bracket operator
\[
\{F_1;F_2\}=\nabla_{\Vec{x}}F_1\cdot\nabla_{\Vec{X}}F_2-\nabla_{\Vec{X}}F_1\cdot\nabla_{\Vec{x}}F_2.
\]
\par

Terms $\widetilde{H}_{0,m}$ in the right side of Eq.~(\ref{homoeq}) are known from previous computations obtained from the general recursion
\begin{equation} \label{triangle}
F_{n,q}=F_{n+1,q-1}+\sum_{0\le{m}\le{n}}{n\choose{m}}\,\{F_{n-m,q-1};W_{m+1}\},
\end{equation}
in which terms $F_{i,j}\equiv{F}_{i,j}(\Vec{x},\Vec{X})$ are replaced by corresponding terms $H_{i,j}$ generated from the known coefficients $H_{m,0}$ of the original Hamiltonian. Terms $H_{0,m}$ in the right side of Eq.~(\ref{homoeq}) are chosen at our will, but with the obvious constraint that the homological equation, which will be, in general, a partial differential equation, can be solved for $W_m$. More properly, $H_{0,m}$ must contain  all the elements of the homological equation that pertain to the kernel of the Lie derivative so that the right side of Eq.~(\ref{homoeq}) only consists of terms of the image of the Lie derivative. In this way, the solution of $W_m$ from Eq.~(\ref{homoeq}) is guaranteed \citep{MeyerHall1992}.
\par

Finally, the new Hamiltonian is obtained by interchanging old and new variables in the $H_{0,m}$ terms, to immediately get Eq.~(\ref{Hnew}).
\par

Once the generating function is computed, the transformation (\ref{Liet}) is explicitly constructed by replacing scalar terms $F_{i,j}$ in Eq.~(\ref{triangle}) by vectorial terms terms $\Vec{x}_{i,j}$ (resp. $\Vec{X}_{i,j}$), starting from $\Vec{x}_{0,0}=\Vec{x}$, and $\Vec{x}_{m,0}=\Vec{0}$ when $m>0$.
\par

It is worth to remark here the power of Deprit's recursion, Eq.~(\ref{triangle}). Since it generally applies to any function of the canonical variable set $(\Vec{x},\Vec{X})$, it can be applied to obtain the transformation (\ref{Liet}) ---typically from mean to osculating elements--- in a different set of canonical variables, say $\Vec{y}=\Vec{y}(\Vec{x},\Vec{X})$, $\Vec{Y}=\Vec{Y}(\Vec{x},\Vec{X})$. But Eq.~(\ref{triangle}) is not restricted to the canonical realm, and, therefore can also be used for obtaining the transformation (\ref{Liet}) in non-canonical variables, as the classical Keplerian orbital elements or different sets of non-singular variables based on them, cf. \citep{DepritRom1970}. In fact, Deprit's recursion applies to vectorial flows, a fact that allowed different extensions of the Lie transforms method to non-canonical perturbations \citep{Kamel1969,Kamel1970,Kamel1971,Henrard1970Lie} ---yet alternative approaches can, of course, be taken in this case \citep{CellettiLhotka2012,LhotkaCelletti2013}. Furthermore, since Poisson brackets are invariant with respect to canonical transformations, they can be evaluated in the more convenient set of canonical variables without constraint to the working set $(\Vec{x},\Vec{X})$.
\par

\section{Brouwer's short-period averaging} \label{s:Brouwer}

When the Lie transformation (\ref{Liet}) of the main problem Hamiltonian (\ref{Hammain}) is approached in Delaunay variables, the zeroth order term of Eq.~(\ref{Hold}) is the Keplerian (\ref{keplerian})
\[
H_{0,0}=-\frac{\mu^2}{2L^2},
\] 
and the Lie derivative (\ref{LieDerivative}) takes the extremely simple form
\begin{equation}
\mathcal{L}_0(W_m)\equiv{n}\,\frac{\partial{W}_m}{\partial\ell},
\end{equation}
where $n=L/a^2$ is the mean motion. The homological equation (\ref{homoeq}) is then solved by the simple quadrature
\begin{equation} \label{sohomeq}
W_m=\frac{1}{n}\int(\widetilde{H}_{0,m}-H_{0,m})\,\mathrm{d}\ell,
\end{equation}
which is undetermined up to an arbitrary function $A_m$ such that $\partial{A}_m/\partial\ell=0$.
\par

At the first order, making $n=0$ and $q=1$ in Eq.~(\ref{triangle}), Deprit's recursion yields
\begin{equation} \label{tl1}
H_{0,1}=H_{1,0}+\{H_{0,0};W_{1}\},
\end{equation}
from which, by comparison to Eq.~(\ref{homoeq}), $\widetilde{H}_{0,1}=H_{1,0}$. Hence,
\begin{equation} \label{kt1}
\widetilde{H}_{0,1}=\frac{\mu}{r}\frac{1}{2}C_{2,0}\frac{R_\oplus^2}{r^2}\left[1-\frac{3}{2}s^2+\frac{3}{2}s^2\cos(2f+2\omega)\right],
\end{equation}
and different choices of $H_{0,1}$ can be made while finding $W_1$ from Eq.~(\ref{sohomeq}). 
\par

\citeauthor{Brouwer1959}'s (\citeyear{Brouwer1959}) plan is to remove first the short-period effects from the Hamiltonian. Therefore, he selects $H_{0,1}$ in such a way that it cancels the terms of $\widetilde{H}_{0,1}$ that are not periodic in the mean anomaly. To this effect, $\widetilde{H}_{0,1}$ is averaged over the mean anomaly, and, to avoid expansions of the true anomaly in Eq.~(\ref{kt1}) in powers of the eccentricity, the differential relation
\begin{equation}
a^2\eta\,\mathrm{d}\ell=r^2\,\mathrm{d}f,
\end{equation}
which is derived from the preservation of the angular momentum in the Kepler problem, is used to obtain
\begin{equation} \label{H01Brouwer}
H_{0,1}=\frac{1}{2\pi}\int_0^{2\pi}\widetilde{H}_{0,1}\frac{r^2}{a^2\eta}\,\mathrm{d}f
=\frac{\mu}{a}\frac{1}{4}C_{2,0}\frac{R_\oplus^2}{p^2}\eta(2-3s^2).
\end{equation}

Then, the first term of the homological equation is computed from Eq.~(\ref{sohomeq}) also in closed form with respect to the eccentricity
\[
W_1=\frac{1}{n}\left(-H_{0,1}\ell+\int\widetilde{H}_{0,1}\frac{r^2}{a^2\eta}\,\mathrm{d}f\right),
\]
yielding
\begin{eqnarray} \label{W1Brouwer}
W_1 &=& G\frac{1}{8}C_{2,0}\frac{R_\oplus^2}{p^2} \left[\left(4-6 s^2\right) (\phi+e \sin{f} ) \right.
\\&& \nonumber \left. +3es^2\sin(f+2\omega)+3s^2\sin(2f+2\omega)+es^2\sin(3f+2\omega)\right] + A_1,
\end{eqnarray}
where
\begin{equation} \label{center}
\phi=f-\ell,
\end{equation}
is the equation of the center, a $2\pi$ periodic, yet not trigonometric function of the mean anomaly, whose closed form treatment in the computation of higher orders of Brouwer's original solution led to the invention of Hamiltonian simplification theories.
\par

Brouwer's choice $A_1=0$ for the arbitrary ``constant'', means that $W_1$ comprises both short- and long-period effects ---as noted by \cite{Kozai1962} after averaging Eq.~(\ref{W1Brouwer}) over the mean anomaly. But this fact is irrelevant to Brouwer's aim of computing the secular terms of the solution by applying a further transformation that removes the long-period terms from the mean elements Hamiltonian.
\par

The first order terms of the mean to osculating transformation (the short-period corrections) are now easily derived from Eq.~(\ref{W1Brouwer}), among which attention is only paid here to the short-period corrections to the inclination $I_{0,1}=\{I,W_1\}$, as obtained from Eq.~(\ref{triangle}). Thus, recalling that $H$ is an integral of the main problem, $I=I(G)$ and the Poisson bracket is easily evaluated to yield
\begin{equation} \label{i01Brouwer}
I_{0,1}=-\frac{1}{4}C_{2,0}\frac{R_\oplus^2}{p^2}\left[3e\cos(f+2\omega)+3\cos(2f+2\omega)+e\cos(3f+2\omega)\right],
\end{equation}
where the true anomaly and the elements in the right side of Eq.~(\ref{i01Brouwer}), viz. $p$, $e$, and $\omega$, are functions of the Delaunay prime variables.
\par

At the second order, the known terms of the homological equation are computed using, again, Deprit's recursion (\ref{triangle}), to get
\begin{equation} \label{tl2}
\widetilde{H}_{0,2}=H_{2,0}+\{H_{0,1};W_{1}\}+\{H_{1,0};W_{1}\},
\end{equation}
where $H_{2,0}\equiv0$ in view of the perturbation arrangement of the main problem. The appearance of the equation of the center in the first order term of the generating function, Eq.~(\ref{W1Brouwer}), results in the consequent generation of periodic terms in the right side of the second order homological equation that cannot be reduced to a finite series of sine and cosine functions. This is, in particular, the case of the term
\begin{equation} \label{chi}
\chi=\frac{9}{8}\frac{\mu}{r}C_{2,0}^2\frac{R_\oplus^4}{p^4}\frac{p^2}{r^2}s^2\left(4-5s^2\right)\phi\sin(2f+2\omega),
\end{equation}
which stems from the evaluation of the Poisson bracket $\{H_{1,0};W_{1}\}$ in Eq.~(\ref{tl2}).
\par

The closed form averaging of terms of the form of Eq.~(\ref{chi}) over the mean anomaly is anyway feasible \citep{Metris1991}, and Brouwer himself managed to compute the second order contribution
\[
\langle\chi\rangle_\ell=\frac{3}{16}\frac{\mu}{p}C_{2,0}^2\frac{R_\oplus^4}{p^4}
\frac{1-\eta}{1+\eta}(1+2\eta)\eta^3s^2\left(4-5s^2\right)\cos2\omega,
\]
of Eq.~(\ref{chi}) to the mean elements Hamiltonian. Still, he didn't pursue the indefinite integration of these kind of terms in closed form, which would allow the solution of the homological equation for the second order term of the generating function, and, in consequence, for the computation of second order short-period corrections \citep{Kozai1962}.
\par

\section{A simplification that looks like a complication} \label{s:simcom}

On the other hand, \cite{Deprit1981} noticed that different solutions of the homological equation can be found by relaxing the strong requirement of removing all the short-period effects from the first order Hamiltonian. Indeed, he used the conic relation (\ref{conic}) to rewrite Eq.~(\ref{kt1}) in the form
\begin{eqnarray} \label{kt1D}
\widetilde{H}_{0,1} &=& \frac{\mu}{p}\frac{1}{8}C_{2,0}\frac{R_\oplus^2}{r^2}
\left[ \left(4-6 s^2\right) (1+e\cos{f}) \right.
\\&& \nonumber \left. +3s^2e\cos(f+2\omega)+6s^2\cos (2 f+2\omega)+s^2e \cos (3 f+2\omega) \right].
\end{eqnarray}
If now $H_{0,1}$ is chosen by taking the terms of $\widetilde{H}_{0,1}$ which are free from the explicit appearance of $f$, viz.
\begin{equation} \label{H01Deprit}
H_{0,1}=\frac{\mu}{p}\frac{1}{4}C_{2,0}\frac{R_\oplus^2}{r^2}\left(2-3s^2\right),
\end{equation}
the first order Hamiltonian still bears short-period effects in the prime variables, but has been \emph{simplified} by, using Deprit's lingo, removing the parallactic factors $1/r^m$ with $m>2$.
\par

The choice of Eq.~(\ref{H01Deprit}) for the new Hamiltonian term instead of Brouwer's selection in Eq.~(\ref{H01Brouwer}), cancels the appearance of the radius in the integrand of the first order homological equation when the true anomaly is taken as the independent variable, in this way trivializing the obtention of the closed form solution. That is, Deprit's selection of $H_{0,1}$ in the form of Eq.~(\ref{H01Deprit}) still keeps the terms of the disturbing function pertaining to the kernel of the Lie derivative, in addition to carrying some terms of its image. The integration of the homological equation is then achieved by simply changing cosines by sines, taking into account the rational factors that may stem from the multiplicity of the argument $f$ in the trigonometric function. Hence,
\begin{eqnarray} \label{W1Deprit}
W_1 &=& G\frac{1}{8}C_{2,0}\frac{R_\oplus^2}{p^2} \left[\left(4-6 s^2\right)e \sin{f} \right.
\\&& \nonumber \left. +3es^2\sin(f+2\omega)+3s^2\sin(2f+2\omega)+es^2\sin(3f+2\omega)\right] + A_1,
\end{eqnarray}
where one may notice that the equation of the center (\ref{center}) is now absent. Because of that, Deprit's reduction of parallactic terms is straightforwardly extended to any order, which he rigorously proved by furnishing the \emph{elimination of the parallax} simplification with the necessary algebraic foundations \citep{Deprit1981,DepritMiller1989,DepritFerrer1989}.
\par

In particular, elementary computations yield the second order Hamiltonian term
\begin{eqnarray} \label{H02Deprit}
H_{0,2} &=& \frac{1}{64}C_{2,0}^2\frac{R_\oplus^4}{p^4}\frac{p}{r}\frac{\mu}{r}\left[
-80+168s^2-84s^4 \right.
\\&& \nonumber \left. -3\left(8-8 s^2-5 s^4\right)e^2 +6 e^2 \left(14-15 s^2\right) s^2 \cos2\omega \right].
\end{eqnarray}
Again, the factor $1/r^2$ in Eq.~(\ref{H02Deprit}) makes the computation of the second order term of the generating function, $W_2$, in closed form to be as simple as before when the homological equation is integrated in the true anomaly, from which the short-period corrections are obtained up to the second order. Recall that $p$, $r$, $s$, $e$, and $\omega$ in Eq.~(\ref{H02Deprit}) are functions of the Delaunay prime variables.
\par

While Deprit's elimination of the parallax does not remove all the short-period effects from the new Hamiltonian, it eases notably this task in a next Delaunay normalization, which is a Lie transformation that converts $L$ into a formal integral of the transformed Hamiltonian in new variables \citep{Deprit1982}. When it is carried out, one finds that the coupling of the equation of the center with standard trigonometric functions of the true anomaly is delayed to the third order, and, therefore, the second order corrections of the short-period elimination missed in Brouwer's solution are now more easily obtained.
\par

The simplifications operated by the elimination of the parallax in the solution of the main problem by preventing the appearance of the equation of the center in the first order term of the generating function are commonly attributed to the elimination of the parallactic terms from the disturbing function, as the name of the transformation induces to think. But they are not. On the contrary, the simplifications are derived from the fact that the elimination of the parallax produces in turn the elimination of the argument of the latitude from the first order Hamiltonian, a fact that is better noticed when the functions that compose Eq.~(\ref{H01Deprit}) are viewed through the filter of polar variables. When this is done, $r$ is directly a canonical variable, as well as it happens to the argument of the latitude $\theta$. The fact that $\theta$ is a cyclic variable in Eq.~(\ref{H01Deprit}) implies that its conjugate couple, the total angular momentum $\Theta$, has been converted into a formal integral of the new Hamiltonian in primed variables. Therefore, up to first order effects, the main problem after the elimination of the parallax evolves with constant inclination.
\par

Incidentally, the reader may be interested in checking that the short-period corrections to the inclination derived from Eq.~(\ref{W1Deprit}) are exactly the same as those given in Eq.~(\ref{i01Brouwer}). That is, up to first order effects, the elimination of the parallax removes all the short-period effects in what respects to the orbital inclination ---which is not at all the case of the other orbital elements.
\par

But the removal of the argument of the latitude is achieved by \emph{all} radial intermediaries, so the question of whether the elimination of the parallax should be preferred to other simplification alternatives naturally emerges. To further illustrate this issue, a new Lie transformation is proposed which increases the exponent of the parallactic factors, as opposed to their elimination. It will be shown that the new transformation results in a simplification, but also a possible complication. More precisely, it produces the same simplifications as the elimination of the parallax in what respects to coping with the equation of the center, but this pertinent simplification may be at the cost of introducing additional short-period perturbations in the orbital elements apart from the inclination. The latter, however, only involves the treatment of longer Poisson series in the consequent Delaunay normalization that removes completely the short-period terms. Because expansions in the eccentricity are not involved, the number of additional periodic terms is small and the task is not viewed as a major concern since the appearance of the so-called Poisson series processors \citep{DanbyDepritRom1965,Rom1970,Giorgilli1979,Henrard1988}. Later on, in Section \ref{s:Cid}, it will be shown that a particular case of this new application of the Lie transforms method ---which is, in fact, an instance of the multiple representation of functions that had already been identified in reference to the elimination of the latitude simplification \citep{CoppolaPalacian1994}--- provides a true simplification in the process of removing the short-period terms.
\par

Let us come back to Deprit's choice in Eq.~(\ref{H01Deprit}), and turn to the conic equation, as he did, but now to decompose the square of the inverse of the radius in the form
\begin{equation} \label{divertimento}
\frac{1}{r^2}=\frac{1}{r^2}\frac{p^2}{r^2}\frac{2}{2+e^2}
-\frac{1}{r^2}\frac{2}{2+e^2}\left(2e\cos{f}+\frac{1}{2}e^2\cos2f\right).
\end{equation}
Clearly, since the second term of Eq.~(\ref{divertimento}) is factored by $1/r^2$ and is expressed as a Fourier series of the true anomaly, it does not cause any trouble if left out of the new Hamiltonian, to be incorporated, after integration in the true anomaly, to the generating function. That is the second term of Eq.~(\ref{divertimento}) pertains to the image of the Lie derivative while the first term still contains all the terms pertaining to the kernel.
\par

Therefore, the new Hamiltonian term
\begin{equation} \label{H01Lara}
H_{0,1}=\frac{1}{2}C_{2,0}\frac{R_\oplus^2}{p^2}\frac{\mu}{r}\frac{p^3}{r^3}\frac{1}{2+e^2}(2-3s^2),
\end{equation}
in which the inverse of the radius is raised to the fourth power, is chosen now at the first order. As far as the argument of the latitude is also an ignorable variable of Eq.~(\ref{H01Lara}), the simplification objective in what respects to subsequent manipulations involving the equation of the center is equally achieved. This fact is checked by computing the first order term of the new generating function by solving, once more, Eq.~(\ref{sohomeq}).
\par

Thus,
\begin{eqnarray} \label{W1Lara}
W_1 &=& G\frac{1}{8}C_{2,0}\frac{R_\oplus^2}{p^2} \left[
\frac{e^2-2}{e^2+2}\left(4-6 s^2\right)e\sin{f}
-\frac{1}{2+e^2}\left(2-3s^2\right)e\sin2f \right. \qquad
\\&& \nonumber \left.+3es^2\sin(f+2\omega)+3s^2\sin(2f+2\omega)+es^2\sin(3f+2\omega)\right] + A_1,
\end{eqnarray}
which has one more trigonometric term ---the term in $\sin2f$--- than its counterpart in Eq.~(\ref{W1Deprit}). Therefore, longer periodic corrections are expected now to be obtained than those provided by the elimination of the parallax. And this is exactly the case, except for the short-period corrections to the inclination, which remain the same as those given in Eq.~(\ref{i01Brouwer}). So, while the new choice of the Hamiltonian produces a simplification by avoiding the appearance of the equation of the center in the generating function, it also complicates the short-period corrections by the introduction of new trigonometric terms which are related to the increased power of the inverse of the radius of the new selection made in Eq.~(\ref{H01Lara}) for $H_{0,1}$.
\par

Progressing to higher orders with the requirement that the new Hamiltonian terms keep $1/r^4$ as a factor, instead of $1/r^2$, is straightforward using the decomposition in Eq.~(\ref{divertimento}). After reaching the desired order of the new Lie transforms simplification, the subsequent removal of short-period terms by means of the standard Delaunay normalization does not involve any additional difficulty with respect to the case in which the preliminary simplification is made by the elimination of the parallax, except form the already mentioned longer trigonometric series involved in the computations.
\par

\section{The neutral case: a new radial intermediary} \label{s:Cid}

Raising the power of the inverse of the radius is, obviously, not constrained to the exponent 4 provided by the decomposition given in Eq.~(\ref{divertimento}), and can be analogously achieved for any other integer exponent. However, in view of increasing the power of the inverse of the radius to higher values produces a corresponding  increment of the number of trigonometric terms of the generating function, keeping this number small seems a reasonable strategy. Fixing this exponent to the value 2 is the successful case of the elimination of the parallax. {When completed at the first order}, the elimination of the parallax provides a useful \emph{natural} intermediary \citep{Deprit1981}. However, when viewed as a \emph{common} intermediary, it preserves less features of the original mean problem Hamiltonian than other intermediaries in the literature. In particular, as pointed out by \cite{DepritFerrer1989}, the radial intermediary of \cite{CidLahulla1971} only departs from the main problem Hamiltonian in those short-period effects contributed by the argument of the latitude. Later, \cite{CidFerrerSeinEchaluce1986} proposed a higher order intermediary that combines the merits of Cid's first order intermediary with additional higher order effects, which they claim to be more representative of the main problem dynamics than \citeauthor{AlfriendCoffey1984}'s (\citeyear{AlfriendCoffey1984}) radial intermediary in which the elimination of the argument of the latitude is obtained after the preliminary elimination of the parallax.
\par

Because it happens that Cid's radial, first order intermediary is a particular instance of the approach discussed in Section \ref{s:simcom} when the exponent of the inverse of the radius is fixed to the value 3, which is precisely the value taken by $1/r$ in the disturbing function of the original main problem Hamiltonian in Eq.~(\ref{Hammain}), this appealing option is further explored in what follows. This case is called neutral because it neither rises nor reduces the exponent of the radius. 
\par

Thus, the square of the inverse of the radius is now decomposed in the form
\begin{equation} \label{diverCid}
\frac{1}{r^2}=\frac{1}{r^2}\frac{p}{r}-\frac{1}{r^2}e\cos{f}.
\end{equation}
Then, Deprit's selection for $H_{0,1}$ in Eq.~(\ref{H01Deprit}) is replaced by
\begin{equation} \label{H01Cid}
H_{0,1}=\frac{1}{4}C_{2,0}\frac{R_\oplus^2}{p^2}\frac{\mu}{r}\frac{p^2}{r^2}(2-3s^2).
\end{equation}
With this new choice and using the $\widetilde{H}_{0,1}$ previously computed in Eq.~(\ref{kt1D}), the homological equation (\ref{sohomeq}) is solved to give
\begin{equation} \label{W1Cid}
W_1=G\frac{1}{8}C_{2,0}\frac{R_\oplus^2}{p^2}s^2\left[
3e\sin(f+2\omega)+3\sin(2f+2\omega)+e\sin(3f+2\omega)\right],
\end{equation}
that is the same as Eq.~(7) of \citep{CidFerrerSeinEchaluce1986} once released from the corset imposed by the algebra of the parallactic state functions used by these authors. The arbitrary constant $A_1$ has been chosen to vanish for simplicity. However, if having a generating function that is composed only of short-period terms were a requirement, the choice
\[
A_1=-G\frac{1}{8}C_{2,0}\frac{R_\oplus^2}{p^2}\frac{1-\eta}{1+\eta}(1+2\eta)s^2\sin2\omega,
\]
would fix the issue, cf.~\citep{LaraSanJuanFolcikCefola2011,LaraSanJuanLopezOchoa2013b}.
\par

At the second order, the known terms that enter the homological equation are given by Eq.~(\ref{tl2}). After evaluation of the involved Poisson brackets using Eqs.~(\ref{H01Cid}) and (\ref{W1Cid}), it is obtained
\begin{eqnarray} \label{Ht02Cid}
\widetilde{H}_{0,2} &=& 
\frac{1}{128}C_{2,0}^2\frac{R_\oplus^4}{p^4}\frac{\mu}{r}\frac{p}{r}\Big\{
6 s^2 \left[e^2 \left(23s^2-16\right)+8 \left(4 s^2-3\right)\right]
\\&& \nonumber
+12 e \left(39 s^2-28\right) s^2 \cos{f}
+6 e^2 \left(23 s^2-16\right) s^2 \cos2f
\\&& \nonumber
+12e^2 \left(14-15 s^2\right)s^2\cos2\omega
+48 e \left(11-12 s^2\right) s^2 \cos (f+2 \omega)
\\&& \nonumber
-24 s^2 \left[e^2 \left(11s^2-10\right)+2 \left(9 s^2-8\right)\right] \cos (2 f+2 \omega)
\\&& \nonumber
-48 e \left(8s^2-7\right) s^2 \cos (3 f+2 \omega)
+12 e^2 \left(6-7s^2\right) s^2 \cos (4 f+2 \omega)
\\&& \nonumber
-15 e^2 s^4 \cos (2 f+4 \omega)
-18 e s^4 \cos (3 f+4 \omega)
-6 \left(e^2-4\right)s^4
\\&& \nonumber
\times\cos (4 f+4 \omega)
+30 e s^4\cos (5 f+4 \omega)
+9 e^2 s^4 \cos (6 f+4 \omega)
\Big\},
\end{eqnarray}
from which the terms
\begin{equation} \label{Qu}
Q=-\frac{3}{64}C_{2,0}^2\frac{R_\oplus^4}{p^4}\frac{\mu}{r}\frac{p}{r}s^2\left[8(3-4s^2)+ (16-23s^2)e^2-2(14-15s^2)e^2\cos2\omega\right],
\end{equation}
are identified as those to be necessarily eliminated previously to the computation of the generating function. That is, all the terms of of Eq.~(\ref{Ht02Cid}) that pertain to the kernel of the Lie derivative are present in Eq.~(\ref{Qu}). Instead of incorporating directly these terms to the new, second order Hamiltonian $H_{0,2}$, $Q$ is fattened up with additional periodic terms of $\widetilde{H}_{0,2}$ using Eq.~(\ref{diverCid}). After that, the new Hamiltonian
\begin{equation} \label{H02Cid}
H_{0,2}=\frac{p}{r}Q(r,\theta,R,\Theta;N),
\end{equation}
is chosen.
\par

The next step is to solve the homological equation to get
\begin{eqnarray} \label{W2Cid}
W_{2} &=& 
\frac{G}{256}C_{2,0}^2\frac{R_\oplus^4}{p^4}\Big\{
-12 e \left(2-e^2\right) \left(16-23 s^2\right) s^2 \sin{f}
\\&& \nonumber
-6 e^2 \left(16-23 s^2\right) s^2 \sin{2f}
+12 e^3 \left(15 s^2-14\right) s^2 \sin (f-2 \omega )
\\&& \nonumber
+12 e s^2 \left[e^2 \left(15 s^2-14\right)-96 s^2+88\right] \sin (f+2\omega )
\\&& \nonumber
-24 s^2 \left[e^2 \left(11 s^2-10\right)+2 \left(9s^2-8\right)\right] \sin (2 f+2 \omega )
\\&& \nonumber
-32 e \left(8 s^2-7\right) s^2 \sin (3 f+2 \omega )
+6 e^2 \left(6-7 s^2\right) s^2 \sin(4 f+2 \omega )
\\&& \nonumber
-15 e^2 s^4 \sin (2 f+4\omega )
-12 e s^4 \sin (3 f+4 \omega )
-3 \left(e^2-4\right) s^4 
\\&& \nonumber
\times\sin (4 f+4 \omega )
+12 e s^4 \sin (5 f+4 \omega )
+3 e^2 s^4 \sin (6 f+4\omega )
\Big\},
\end{eqnarray}
which has the only extra trigonometric term whose argument is $(f-2\omega)$ with respect to the corresponding term of the generating function of the elimination of the parallax. 

Proceeding analogously, the new Hamiltonian term
\begin{eqnarray} \label{H03Cid}
H_{0,3} &=& \frac{9}{1024}C_{2,0}^3\frac{R_\oplus^6}{p^6}\frac{\mu}{r}\frac{p^2}{r^2}s^2
\left\{16 (31 - 26 s^2) s^2 \right.
\\&& \nonumber  - 4 (440 - 1614 s^2 + 1217 s^4) e^2 - 9 (2 - 3 s^2) (16 - 23 s^2) e^4
\\&& \nonumber \left.  +  2 \left[1984 - 5264 s^2 + 3405 s^4 + 
     9 (28 - 72 s^2 + 45 s^4) e^2\right] e^2\cos2\omega\right\},
\end{eqnarray}
is obtained with the help of Eq.~(\ref{diverCid}). As required, it is factored by the third power of the inverse of the radius.
\par

Finally, following \citeauthor{Izsak1963AJ}'s (\citeyear{Izsak1963AJ}) advice, in order to get a compact form of the periodic corrections, the terms $W_1$ and $W_2$ of the generating function are reformulated in polar variables. Namely, Eqs.~(\ref{W1Cid}) and (\ref{W2Cid}) are replaced by
\begin{eqnarray}
W_1 &=& G\frac{1}{8}C_{2,0}\frac{R_\oplus^2}{p^2}s^2\left[(3+4\kappa)\sin2\theta-2\sigma\cos2\theta\right],
\\[1ex]
W_2 &=& G\frac{1}{64}C_{2,0}^2\frac{R_\oplus^4}{p^4}\Big\{
3\left(16-23s^2\right)s^2\left(\kappa^2+\sigma^2-\kappa-2\right)\sigma
\\&& \nonumber
-\left[16(13-14s^2)-3(6-7s^2)\kappa+6(14-15s^2)(\kappa^2-\sigma^2)\right]s^2\sigma\cos2\theta
\\&& \nonumber
+3s^4(2+3\kappa)\sigma\cos4\theta
\\&& \nonumber
+2\left[6(8-9s^2)+16(10-11s^2)\kappa-6(14-15s^2)\kappa\sigma^2+\frac{3}{4}(46-51s^2)\kappa^2 \right.
\\&& \nonumber \left.
+\frac{3}{4}(34-37s^2)\sigma^2\right]s^2\sin2\theta
+\frac{3}{4}\left(4-5\kappa^2+3\sigma^2\right)s^4\sin4\theta
\Big\},
\end{eqnarray}
in which the projections of the eccentricity vector in the orbital frame, $\kappa=e\cos{f}$ and $\sigma=e\sin{f}$, are substituded by the functions
\begin{equation} \label{kasi}
\kappa=\frac{p}{r}-1, \qquad \sigma=\frac{pR}{\Theta},
\end{equation}
of the polar variables.
\par

Note that, for the lower eccentricity orbits, terms factored by the different powers of the eccentricity in Eqs.~(\ref{H02Cid}) and (\ref{H03Cid}) can be neglected. When this is done, the disturbing function no longer depends on either the eccentricity or the argument of the perigee, which in turn means that, on the one hand, the argument of the latitude has been eliminated, and, on the other, the disturbing function has been released from its dependence on the radial velocity, therefore yielding a radial, higher order intermediary that is formally analogous to Cid's first order intermediary, thus admitting an equally simple integration, cf.~\citep{Ferrandiz1986,LaraGurfil2012}.
\par

In the general case, however, the eccentricity cannot be neglected. Then, to complete the computation of the radial intermediary, the dependence on the argument of the latitude must be removed from the second and third order terms of the new Hamiltonian. Because the appearance of the argument of the latitude in Eqs~(\ref{H02Cid}) and (\ref{H03Cid}) is only due to the term
\begin{equation} \label{cos2w}
e^2\cos2\omega=(\kappa^2-\sigma^2)\cos2\theta+2\kappa\sigma\sin2\theta,
\end{equation}
the (partial) normalization of the simplified Hamiltonian can be done analogously to \citeauthor{AlfriendCoffey1984}'s (\citeyear{AlfriendCoffey1984}) elimination the perigee, yet keeping factors $1/r^3$ in the new Hamiltonian, contrary to $1/r^2$. Besides, the normalization procedure is more conveniently approached in Delaunay variables, cf.~\citep{LaraSanJuanLopezOchoa2013c}.
\par

Thus, a new canonical transformation 
\begin{equation} \label{Liet2}
\mathcal{T}_1:(\Vec{x}',\Vec{X}')\longrightarrow(\Vec{x}'',\Vec{X}'';\epsilon),
\end{equation}
from the primed Delaunay variables in which the Hamiltonian
\[
\mathcal{K}=\sum_{m=0}^3\frac{\epsilon^m}{m!}K_{m,0},
\]
with $K_{0,0}=H_{0,0}$, $K_{1,0}=H_{0,1}$, as given in Eq.~(\ref{H01Cid}), $K_{2,0}=H_{0,2}$, as given in Eq.~(\ref{H02Cid}), and $K_{3,0}=H_{0,3}$, as given in Eq.~(\ref{H03Cid}), is expressed, to new, double prime variables
in which the transformed Hamiltonian
\begin{equation} \label{Ksin}
\mathcal{Q}=\mathcal{K}(\Vec{x}'(\Vec{x}'',\Vec{X}'';\epsilon),\Vec{X}'(\Vec{x}'',\Vec{X}'';\epsilon);\epsilon)\equiv\sum_{m=0}^3\frac{\epsilon^m}{m!}K_{0,m}(\Vec{x}'',\Vec{X}''),
\end{equation}
is free from the argument of the latitude, is computed with the Lie transforms method as follows.
\par

The first order term of the new Hamiltonian remains unchanged in view of it is already free form the argument of the perigee. Hence,
\begin{equation} \label{K01Cid}
K_{0,1}=K_{1,0}.
\end{equation}
Then, the right side of the homological equation (\ref{homoeq}) identically vanishes and, in consequence, Eq.~(\ref{sohomeq}) yields
\begin{equation} \label{W1perigee}
W_1=A_1(-,g,-,L,G,H),
\end{equation}
which, by now, is left undetermined.
\par

At the second order, the terms $\widetilde{K}_{0,2}$ involved in the solution of the homological equation must be computed from the analog to Eq.~(\ref{tl2}), that is
\begin{equation} \label{tl2bis}
\widetilde{K}_{0,2}=K_{2,0}+\{K_{0,1};W_{1}\}+\{K_{1,0};W_{1}\},
\end{equation}
but now, since $W_{1}$ is yet undetermined, its partial derivatives will show up in the evaluation of the Poisson brackets. After the usual preprocessing, it is obtained
\begin{eqnarray} \label{Kt02Cid}
\widetilde{K}_{0,2} &=& 
-\frac{3}{64}C_{2,0}^2\frac{R_\oplus^4}{p^4}\frac{\mu}{r}\frac{p}{r}s^2\Big\{
\left[8 \left(3-4 s^2\right)+\left(16-23 s^2\right)e^2\right] (1+e\cos{f})
\\&& \nonumber -\left(14-15 s^2\right)e^2\left[e \cos(f-2\omega)+2 \cos2\omega+e \cos (f+2\omega)\right]
\Big\}
\\&& \nonumber
-\frac{3}{8\eta^3}C_{2,0}\frac{R_\oplus^2}{p^2}\frac{\mu}{r}\frac{p}{r}(2-3s^2)e\left[(4+e^2)\sin{f}+4e\sin2f+e^2\sin3f\right]\frac{\partial{A}_1}{\partial{L}}
\\&& \nonumber +\frac{3}{2}C_{2,0}\frac{R_\oplus^2}{p^2}\frac{G}{r^2}\Big\{
(4-5s^2)+\frac{1}{4e}\left[8-12 s^2+(14-17 s^2)e^2\right]\cos{f}
\\&& \nonumber+(2-3s^2)\cos2f+\frac{1}{4}(2-3s^2e)\cos3f
\Big\}\frac{\partial{A}_1}{\partial{g}},
\end{eqnarray}
from which it is easy to identify the terms that could introduce secular effects in the computation of the second order term of the generating function, as did happen to \cite{CidLahulla1969,CidLahulla1971}, like those summands that do not depend explicitly on the true anomaly.
\par

Recalling that the partial derivatives of $A_1$ do not involve the mean anomaly, these terms can be obtained from simple inspection of Eq.~(\ref{Kt02Cid}). They are
\begin{eqnarray*}
&& 
-\frac{3}{64}C_{2,0}^2\frac{R_\oplus^4}{p^4}\frac{\mu}{r}\frac{p}{r}s^2
\left[8 \left(3-4 s^2\right)+\left(16-23 s^2\right)e^2\right]
\\
&& 
-\frac{3}{32}C_{2,0}^2\frac{R_\oplus^4}{p^4}\frac{G^2}{r^2}s^2(14 - 15 s^2) e^2\cos2\omega
+\frac{3}{2}\frac{G}{r^2}C_{2,0}\frac{R_\oplus^2}{p^2}(4-5s^2)\frac{\partial{A}_1}{\partial{g}},
\end{eqnarray*}
where it is evident that the term in the first row can be incorporated to the new Hamiltonian for it does not depend either on the argument of the perigee or explicitly on the true anomaly. Following \citeauthor{AlfriendCoffey1984} (\citeyear{AlfriendCoffey1984}), the introduction of secular terms in the generating function is avoided by canceling the second row. That is, by solving $A_1$ from the differential equation in Delaunay variables
\begin{equation} \label{A1Cid}
(4-5s^2)\frac{\partial{A}_1}{\partial{g}}
-\frac{1}{16}C_{2,0}\frac{R_\oplus^2}{p^2}Gs^2(14 - 15 s^2) e^2\cos2\omega=0.
\end{equation}
Hence, up to an integration constant that now might depend only on the Delaunay action variables,
\[
A_1=-\frac{1}{32}C_{2,0}\frac{R_\oplus^2}{p^2}G\frac{14-15s^2}{4-5s^2}s^2e^2\sin2\omega,
\]
which is the same ``constant'' found by \cite{AlfriendCoffey1984}. 
\par

Once the first order term of the generating function $W_1=A_1$ has been determined, and its derivatives have been replaced in Eq.~(\ref{Kt02Cid}), the terms that must be canceled in the solution of the homological equation are fattened up with other periodic terms of the homological equation, as it has been done in the previous simplification, in order to incorporate them to the new Hamiltonian with the desired factor $1/r^3$. After doing that, the second order term of the Hamiltonian with the argument of the latitude removed is
\begin{equation} \label{K02Cid}
K_{0,2}=-\frac{3}{64}C_{2,0}^2\frac{R_\oplus^4}{p^4}\frac{\mu}{r}\frac{p^2}{r^2}s^2
\left[8 \left(3-4 s^2\right)+\left(16-23 s^2\right)e^2\right].
\end{equation}

Next, the second order term of the generating function is computed but, again, only up to an arbitrary function $A_2\equiv{A}_2(-,g,-,L,G,H)$, to yield
\begin{eqnarray} \label{CidW2}
W_2 &=& \frac{1}{128}GC_{2,0}^2\frac{R_\oplus^4}{p^4}es^2\frac{14-15s^2}{4-5s^2}(2-3s^2)\big[3e^2\sin(f-2\omega)
\\&& \nonumber
-e^2\sin(3f+2\omega)-6e\sin(2f+2\omega)-12\sin(f+2\omega)\big]+A_2.
\end{eqnarray}
This arbitrary function is properly determined only at the third order. Proceeding as before, it is obtained
\begin{eqnarray*}
A_2 &=&-\frac{1}{512}C_{2,0}^2\frac{R_\oplus^4}{p^4}\frac{G}{4 - 5 s^2}\left\{
e^2 s^2\Big[4 (824 - 1997 s^2 + 1215 s^4) + \frac{14 - 15 s^2}{4 - 5 s^2} \right.
\\&& \left. \times (56 - 36 s^2 - 45 s^4) e^2\Big]\sin2\omega
- \frac{(14 - 15 s^2)^2}{2(4 - 5 s^2)^2} (13 - 15 s^2) e^4 s^4\sin4\omega
\right\},
\end{eqnarray*}
which makes $W_2$ to become completely determined. Now $A_2$ is different from the one that would be obtained with the original elimination of the perigee \citep{AlfriendCoffey1984}.
\par

The computations are stopped after calculating the third order Hamiltonian
\begin{eqnarray} \label{K03Cid}
K_{0,3} &=&
\frac{9}{1024}C_{2,0}^3\frac{R_\oplus^6}{p^6}\frac{\mu}{r}\frac{p^2}{r^2}\frac{s^2}{4-5s^2}
\left[16 s^2 (4 - 5 s^2)^2 (31 - 26 s^2) \right.
\\&& \nonumber - 2 (4 - 5 s^2) (3520 - 17116 s^2 + 25456 s^4 - 11945 s^6) e^2 
\\&& \nonumber \left.- (3040 - 15116 s^2 + 29176 s^4 - 25815 s^6 +8775 s^8) e^4 \right].
\end{eqnarray}
and the Lie transforms procedure ends by interchanging prime by double primed variables in Eqs.~(\ref{K01Cid}), (\ref{K02Cid}), and (\ref{K03Cid}).
\par

In spite of the use of Delaunay variables made throughout the whole simplification process of the elimination of the perigee, that the argument of the latitude has been removed from the simplified Hamiltonian (\ref{Ksin}) is clearly disclosed when viewed through the filter of polar variables. Indeed, the eccentricity in Eqs.~(\ref{K02Cid}) and (\ref{K03Cid}) is only a function of the radius, the radial velocity, and the total angular momentum. The integration of the new radial intermediary can be approached by different methods \citep[see][for instance]{FrancoPalacios1990,AbadSanJuanGavin2001}.

\section{Conclusions}

In the half a century elapsed since Deprit's approach to Hamiltonian perturbations by Lie transforms appeared in print in the first issue of the Celestial Mechanics journal, a series of successful Hamiltonian and non-Hamiltonian Lie transforms simplifications ---most of which were published also in the same journal--- have proven the power of the method. The new Lie transforms application proposed here, which for sure will not end the series, has been useful in showing that the intricacies introduced by the equation of the center in the removal of short-period effects from a zonal Hamiltonian, are related to the short-period effects that disturb the Keplerian motion in the out of plane direction, instead of its general departure from a quasi-Keplerian system. 
\par

The neutral case of the new transformation, in which neither the parallax is eliminated nor the exponents of the inverse of the radius are raised, led to the discovery of a new radial intermediary for low Earth orbits. The new solution extends the recognized merits of Cid's radial intermediary to higher orders, while enjoying an analogous simple solution.

\section*{Compliance with ethical standards}
\paragraph{Conflict of interest} The author declares that he has no conflict of interest.

\section*{Acknowledgements}

Partial support by the Spanish State Research Agency and the European Regional Development Fund under Projects ESP2016-76585-R and ESP2017-87271-P (AEI/ ERDF, EU) is recognized.

\end{document}